\documentclass[final,leqno,onefignum,oneeqnum]{siamltex704}

\usepackage{amsmath}
\usepackage{amssymb}
\usepackage{amsfonts}
\usepackage{graphicx}
\usepackage{color}
\numberwithin{equation}{section}
\newtheorem{teo}{Theorem}[section]
\newtheorem{lem}[teo]{Lemma}

\newtheorem{oss}[teo]{{\it Remark}}

\newcommand{\ve}{\varepsilon}
\newcommand{\ffi}{\varphi}
\newcommand{\deb}{\rightharpoonup}

\newcommand{\pical}{\mathcal{P}}
\newcommand{\lcal}{\mathcal{L}}

\newcommand{\R}{\mathbb{R}}
\newcommand{\F}{\mathfrak{F}}
\newcommand{\tto}{\rightarrow}

\title{A Mass Transportation Model for the Optimal Planning of an Urban Region\thanks{Received by the editors December 4, 2003; accepted for publication (in revised form) November 29, 2004; published electronically October 27, 2005. This work is part of the European Research Training Network ``Homogenization and Multiple Scales''(HMS{\rm 2000}) under contract HPRN-2000-00109. 
\URL sima/37-2/43831.html}}

\author{Giuseppe Buttazzo\thanks{Dipartimento di Matematica, Universit\`{a} di Pisa, via Buonarroti 2, 56127 Pisa, Italy\break 
(buttazzo@dm.unipi.it). The first author also acknowledges the support of the project ``Problemi di Ottimizzazione in Teoria del Trasporto ed Applicazioni a Problemi di Pianificazione Urbana'' of the Italian GNAMPA and of the project ``Calcolo delle Variazioni'' of the Italian Ministry of Education.} 
\and Filippo Santambrogio\thanks{CEREMADE, Universit\'e Paris-Dauphine, Place de Lattre de Tassigny, 75775 Paris cedex 16, France (filippo@ceremade.dauphine.fr).}}

\begin{document}

\slugger{sima}{2005}{37}{2}{514--530}

\setcounter{page}{514}
\maketitle

\begin{abstract}
We propose a model to describe the optimal distributions of residents and services in a prescribed urban area. 
The cost functional takes into account the transportation costs (according to a Monge--Kantorovich-type criterion) and two additional terms which penalize concentration of residents and dispersion of services. The tools we use are the Monge--Kantorovich mass transportation theory and the theory of nonconvex functionals defined on measures.
\end{abstract}



\begin{DOI}
10.1137/S0036141003438313
\end{DOI}

\pagestyle{myheadings}
\thispagestyle{plain}
\markboth{GIUSEPPE BUTTAZZO AND FILIPPO SANTAMBROGIO}{MASS TRANSPORTATION AND OPTIMAL PLANNING OF URBAN REGIONS}                                                      

\section{Introduction}\label{sec1}
The efficient planning of a city is a tremendously complicated problem, both for the high number of parameters which are involved as well as for the several relations which intervene among them (price of the land, kind of
industries working in the area, quality of the life, prices of transportations, geographical obstacles, etc.). Perhaps a careful description of the real situations could be only obtained through evolution models which take into account the dynamical behavior of the different parameters involved.

Among the several ingredients in the description of a city which are considered by urban planners, two of the most important are the distribution of residents and the distribution of services  (working places, stores, offices, etc.). These two densities have to be treated in a different way due to their features (see facts ii) and iii) below). Several interesting mathematical models for the description of the equilibrium structure of these two elements of a city have been studied in the spatial economical literature (see for instance the classical text by Fujita \cite{Fuj} and the more recent paper by Lucas and Rossi-Hansberg \cite{LucRoH})). The first time where the Monge--Kantorovich theory of optimal transportation appears and plays an important role, is, to the best of our knowledge, in Carlier and Ekeland \cite{carlierekeland,carlierekeland2}.

All these papers are mostly equilibrium-oriented and this is the main difference with the present work, which wants to be focused on the following aspects: distribution of residents and services, optimization of a global criterion, optimal transport theory.

We consider a geographical area as given, and we represent it through a
subset $\Omega$ of $\R^n$ ($n=2$ in the applications to concrete urban
planning problems). We want to study the optimal location in $\Omega$ of
a mass of inhabitants, which we denote by $\mu$, as well as of a mass of
services, which we denote by $\nu$.
We assume that $\mu$ and $\nu$ are probability measures on $\Omega$. This 
means that the total amounts of population and production are 
fixed as problem data, and this is a difference from the model in \cite{carlierekeland}. 
The measures $\mu$ and $\nu$ represent the unknowns of our problem that have to be found 
in such a way that a suitable total cost functional $\F(\mu,\nu)$ is minimized. The definition of
this total cost functional takes into account some criteria we want the two
densities $\mu$ and $\nu$ to satisfy:

\begin{itemize}
\item[(i)] there is a transportation cost for moving from the residential
areas to the services areas;

\item[(ii)] people do not want to live in areas where the density of
population is too high;

\item[(iii)] services need to be concentrated as much as possible in order
to increase efficiency and decrease management costs.

\end{itemize}

Fact (i) will be described through a Monge--Kantorovich mass transportation
model; the transportation cost will indeed be given by using a $p$-Wasserstein
distance ($p\geq1$). We set
\begin{equation}\label{e1wasse}
T_p(\mu,\nu)=W_p^p(\mu,\nu)=\inf_{\gamma }\left(\int_{\Omega\times\Omega}\!\!|x-y|^p\gamma(dx,dy)\right),
\end{equation}
where the infimum is taken over all possible transport plans $\gamma$ between $\mu$ and $\nu$ (i.e., probabilities 
on the product space having $\mu,\nu$ as marginal measures). We refer to \cite{dispensedivillani} for the whole 
theory on mass transportation. When $p=1$ we are in the classical Monge case, and for this particular case we refer 
to \cite{ambpra} and \cite{evansgangbo}.

Fact (ii) will be described by a penalization functional, a kind of total
unhappiness of citizens due to high density of population, obtained by integrating
with respect to the citizens' density their personal unhappiness.

Fact (iii) is modeled by a third term representing costs for managing services once
they are located according to the distribution $\nu$, taking into account that
efficiency depends strongly on how much $\nu$ is concentrated.

The cost functional we will consider is then
\begin{equation}\label{e1costf}
\F^p(\mu,\nu)=T_p(\mu,\nu)+F(\mu)+G(\nu),
\end{equation}
and thus the optimal location of $\mu$ and $\nu$ will be determined by the
minimization problem
\begin{equation}\label{e1minpb}
\min\big\{\F^p(\mu,\nu)\ :\ \mu,\nu\mbox{ probabilities on }\Omega\big\}.
\end{equation}

In this way, our model takes into consideration only the optimization of a
total welfare parameter of the city, disregarding the effects on each single
citizen. In particular, no equilibrium condition is considered. This may appear as a fault in the model, since the personal welfare of the citizens (depending on the population density of their zone and on the cost of moving from home to
services) could be nonconstant. As a consequence, nonstable optimal solutions
may occur, where some citizens would prefer to move elsewhere in the city in
order to get better conditions. However, this is not the case, since our model
also disregards prices of land and houses in the city, since they do not affect
the total wealth of the area. It turns out that by a proper, market-determined
choice of prices, welfare differences could be compensated and equilibrium
recovered. This fact turns out to be a major difference between our model and the model in
\cite{carlierekeland}, both for the importance given in \cite{carlierekeland} to the variable
represented by the price of land and for the fact that Carlier and Ekeland specifically
look for an equilibrium solution instead of an optimal one.

The present paper, after this introduction, contains three sections. Section \ref{sec2} is
devoted to presenting precise choices for the functionals $F$ and $G$ and
justifying them as reasonable choices. In the same section we also give a
simple existence result for an optimal solution $(\mu,\nu)$ as a starting point
for the rest of the paper. In section \ref{sec3} we consider the functional on $\mu$
obtained by keeping the measure $\nu$ as fixed: in this case the functional $G$
does not play any role, and we obtain a convex minimization problem, which is
interesting in itself. We also obtain some necessary optimality conditions in
the very general case where no assumption is taken on the fixed measure $\nu$.
In section \ref{sec4} we apply these results to the case where $G$ is of the particular form presented in section \ref{sec2}, which forces $\nu$ to be atomic (i.e., services are concentrated in countably many points of the city area $\Omega$). In the case where $\Omega$ is bounded, we give a quite precise description of the solution $(\mu,\nu)$, and then we give an existence result also for the case $\Omega=\R^n$. 

Both in the case $\Omega=\R^n$ and $\Omega$ bounded, optimal choices for $\mu$ and $\nu$ are given by the formation of a certain number of subcities, which are circular areas with a pole of services in the center (an atom for the measure $\nu$) around which the population is distributed with a decreasing radial density.

Since we have considered only a very simplified model, our goal is neither 
to suggest a realistic way to design the ideal city nor to describe in a 
variational way the formation of existing cities. Nevertheless, from the analysis of our optimality results (and in particular from the subcities phenomena we referred to), we can infer some conclusions.

\begin{itemize}
\item Our model is not a proper choice to describe the shape of a single existing city, 
since the delocalization of services we find in an optimal solution does not reflect what 
reality suggests (in fact, we find finitely many disjoint, independent subcities with services only in the center).
\item Our model is likely to be more realistic on a larger scale, when $\Omega$ represents 
a large urban area composed of several cities: in this case every atom of the optimal $\nu$ 
stands for the center of one of them and includes a complex system of services, located downtown, 
whose complexity cannot be seen in this scale.
\item In our model the concentrated measure $\nu$ gives a good representation of the areas 
where services are offered to citizens and not of areas where commodities are produced 
(factories), due to the assumption that no land is actually occupied by the service poles (since $\nu$ is atomic).
\item We do not believe that our model may actually be used to plan a future city or 
to improve the efficiency of an existing one, as a consequence of its oversimplified 
nature. However, we do not exclude the possibility of using it in the planning of less 
complex agglomerations, such as tourist villages, university campuses, etc.
\item We conclude by stressing that the same model may be applied as a first 
simplified approach to other kinds of problems, where we have to choose in some 
efficient way the distributions of two different parameters, the first spread 
and the second concentrated, keeping them as close as possible to each other in some 
mass transportation sense.
\end{itemize}

This last issue, the investigation of optimization problems concerning concentration and transport costs of probability measures, has been the main subject of \cite{PhDthesis}, which also refers to the first version of this paper. Two subsequent works included in the same manuscript have explicitly dealt with urban planning problem, with different transport costs (\cite{CarSan}) and with a different service performance evaluation (\cite{T+GV}, see Section 2). For other issues in urban planning, more related to the presence of a transportation system, we cite two different books. The first one, \cite{ButPraSolSte}, summarizes several works by the italian optimal transport community concerning the location of transport networks in a urban area once $\mu$ and $\nu$ are given. The second we want to mention, on the contrary, is a reference book by Beckmann and Puu, \cite{BecPuu}, dealing with several spatial economics problem and approached trough the ``continuous transportation model'', which is much linked to Monge-Kantorovich theory. One of the problem which is addressed is the repartition of land between the residents on the one hand and the road system on the other. For simplicity the distribution of services is identified with that of residents and another interesting variable, affecting transportation costs, is introduced, standing for road width.

\section{The model}\label{sec2}

We now define the three terms appearing in our functional $\F^p$. We must go through 
the definition of $F$ and $G$, since the first term will be a Monge--Kantorovich transport 
cost, as explained in the previous section. 
For the functional $F$ we take
\begin{equation}\label{defiF}
F(\mu)=\begin{cases}\int_{\Omega}f\!\left(u(x)\right)dx&\text{if }
\mu=u\cdot\lcal^n,\,u\in L^1(\Omega),\\+\infty&\text{otherwise,}\end{cases}
\end{equation}
where the integrand $f:[0,+\infty]\tto[0,+\infty]$ is assumed to
be lower semicontinuous and convex, with $f(0)=0$ and superlinear at
infinity, that is,
\begin{equation}\label{superlinearita}\lim_{t\tto+\infty}\frac{f(t)}{t}=+\infty.\end{equation}
In this form we have a local semicontinuous functional on measures. Without loss of generality, 
by subtracting constants to the functional $F$, we can suppose $f'(0)=0$.
Due to the assumption $f(0)=0$, the ratio $f(t)/t$ is an incremental ratio 
of the convex function $f$ and thus it is increasing in $t$. Then, if we write the functional $F$ as
$$\int_{\Omega}\frac{f(u(x))}{u(x)}u(x)\,dx,$$
we can see the quantity $f(u)/u$, which is increasing in $u$, as the unhappiness of a single 
citizen when he lives in a place where the population density is $u$. Integrating it with respect 
to $\mu=u\cdot\lcal^n$ gives a quantity to be seen as the total unhappiness of the population. 

As far as the concentration term $G(\nu)$ is concerned, we set 
\begin{equation}\label{defiG}G(\nu)=\begin{cases}\sum_{i=0}^{\infty}g(a_i)&\text{if
}\nu=\sum_{i=0}^{\infty}a_i\delta_{x_i},\\+\infty&\text{if $\nu$ is not atomic.}\end{cases}\end{equation}
We require the function $g$ to be subadditive, lower semicontinuous, and such that $g(0)=0$ and 
\begin{equation}\label{gizeroinfinita}\lim_{t\tto 0}\frac{g(t)}{t}=+\infty.\end{equation}
Every single term $g(a_i)$ in the sum in (\ref{defiG}) represents the cost of building
and managing a service pole of size $a_i$, located at the point $x_i\in\Omega$.

In our model, as already pointed out, we fix as a datum the total production of services; moreover, 
in each service pole the production is required as a quantity proportionally depending on its size 
(or on the number of inhabitants making use of such a pole). We may define the productivity $P$ of 
a pole of mass (size) $a$ as the ratio between the production and the cost to get such a production. 
Then we have $P(a)=a/g(a)$ and
$$\sum_{i=0}^{\infty}g(a_i)=\sum_{i=0}^{\infty}\frac{a_i}{P(a_i)}.$$
As a consequence of assumption \eqref{gizeroinfinita} we have that the productivity in very small 
service poles is near $0$.

Notice that in the functional $G$ we do not take into account distances between service poles. Nonlocal functionals involving such distances, taking into account possible cooperation and the consequent gain in efficiency, have been considered in subsequent investigations (see \cite{T+GV,PhDthesis}). The results shown in the next section (since they do not depend on the choice of $G$) turned out to be useful in such a setting.

For the problem introduced in (\ref{e1minpb}), existence results are straightforward, especially 
when we use as an environment a compact set $\Omega$. In fact, functionals of the form of 
both $F$ and $G$ have been studied in a general setting by Bouchitt\'{e} and Buttazzo in \cite{bb1}, 
and lower semicontinuity results were proven.
\begin{teo}
Suppose $\Omega$ is compact, $p\geq1$, and $f$ and $g$ satisfy the conditions listed above. 
Then the minimization problem {\rm (\ref{e1minpb})} has at least one solution.
\end{teo}
\unskip

\begin{proof}
By the direct method of calculus of variations, this result is an easy consequence of the weak-* 
compactness of the space $\pical(\Omega)$, the space of probability measures on $\Omega$ when $\Omega$ 
itself is compact, and of the weak-* semicontinuity of the functional $\mathfrak{F}^p$. The second 
and third terms in (\ref{e1costf}) are, in fact, local semicontinuous functionals 
(due to results in \cite{bb1}), while the first term is nothing but a Wasserstein distance raised 
to a certain power. Since it is known that in compact spaces this distance metrizes the weak-* topology, 
$T_p$ is actually continuous.\qquad
\end{proof}

In \cite{miatesi}, where we first presented the model, other existence results were shown. For instance, the case of a noncompact bounded convex set $\Omega\subset\R^n$ was considered. We will not go through this proof here and will discuss just one existence result in a noncompact setting, obtained as a consequence of a proper use of the optimality conditions presented in the next section.

\section{A necessary condition of optimality}\label{sec3}

In this section we find optimality conditions for probability measures on $\Omega$ minimizing the functional
$$\mathfrak{F}^p_{\nu}(\mu)=T_p(\mu,\nu)+F(\mu).$$
It is clear that if $(\mu,\nu)$ is an optimal pair for the whole functional $\mathfrak{F}^p$, 
$\mu$ is a minimizer for
$\mathfrak{F}^p_{\nu}$. The goal of this section is to derive optimality conditions for $\mathfrak{F}^p_{\nu}$, 
for any $\nu$, without any
link to the minimization of $\F^p$. The main part of the section will be devoted to presenting an approach 
obtained by starting with the easier case $p>1$ and $\nu$ ``regular'' in some sense and then recovering the general case
by an approximation argument. The reason for doing so relies on some conditions ensuring uniqueness properties of 
the Kantorovich potential. Similar approximation arguments were also used in \cite{miatesi}: purely atomic probability 
measures (i.e., finite sums of Dirac masses) were considered first, and then, by approximation, the result was extended 
to any measure $\nu$. At the end of the section we also provide a sketch of a different proof, suggested to us by 
an anonymous referee, which is based on some convex analysis tools and strongly uses the convex structure of the problem. 

For simplicity, let us call domains those sets which are the closure of a nonempty connected open subset of $\R^n$ 
with negligible boundary. From now on $\Omega$ will be a bounded domain and its diameter will be denoted by $D$. 
The function $f$ in \eqref{defiF} will be assumed to be strictly convex and $C^1$, and we will denote by $k$ 
the continuous, strictly increasing function $(f')^{-1}$. Strict convexity of $f$ will
ensure uniqueness for the minimizer of $\mathfrak{F}^p_{\nu}$.

\begin{lem}\label{stimafacileconconv}
If $\mu$ is optimal for $\mathfrak{F}^p_{\nu}$, then for any other probability measure $\mu_1$ with density $u_1$
such that
$\mathfrak{F}^p_{\nu}(\mu_1)<+\infty$, the following inequality holds:
$$T_p(\mu_1,\nu)-T_p(\mu,\nu)+\int_{\Omega} f'(u(x))[u_1(x)-u(x)]dx\geq0.$$
\end{lem}
\unskip

\begin{proof}
For any $\ve>0$, due to the convexity of the transport term, it holds that
\begin{multline*}
T_p(\mu,\nu)+F(\mu)\leq T_p(\mu+\ve(\mu_1-\mu))+F(\mu+\ve(\mu_1-\mu),\nu)\\
\leq T_p(\mu,\nu)+\ve(T_p(\mu_1,\nu)-T_p(\mu,\nu))+F(\mu+\ve(\mu_1-\mu)).
\end{multline*}
Therefore the quantity
$T_p(\mu_1,\nu)-T_p(\mu,\nu)+{\ve}^{-1}\left[F(\mu+\ve(\mu_1-\mu))-F(\mu)\right]$
is nonnegative. If we let $\ve\tto 0$, we obtain the thesis if we prove
$$\limsup_{\ve\tto0}\int\frac{f(u+\ve(u_1-u))-f(u)}{\ve}\,d\lcal^n\leq\int
f'(u)(u_1-u)\,d\lcal^n.$$ 
By using the monotonicity of the incremental ratios of convex functions we can see that, for $\ve<1$,
$$\frac{f(u+\ve(u_1-u))-f(u)}{\ve}\leq f(u_1)-f(u).$$
This is sufficient in order to apply Fatou's Lemma, since the quantities $F(\mu)$ and $F(\mu_1)$ are finite.\qquad
\end{proof}

\begin{lem}\label{muriempietuttosenilofa}
Let us suppose $\nu=\nu^s+v\cdot\lcal^n$, with $v\in L^{\infty}(\Omega),\,\nu^s\bot\lcal^n,\,v>0$ a.e. in $\Omega$. 
If $\mu$ is optimal for $\mathfrak{F}^p_{\nu}$, then $u>0$ a.e. in $\Omega$.
\end{lem}
\unskip

\begin{proof}
The lemma will be proven by contradiction. We will find, if the set $A=\left\{u=0\right\}$ is not negligible, 
a measure $\mu_1$ for which Lemma \ref{stimafacileconconv} is not verified. Let $N$ be a Lebesgue-negligible 
set where $\nu^s$ is concentrated and $t$ is an optimal transport map between $\mu$ and $\nu$. Such an optimal 
transport exists, since $\mu\ll\lcal^n$. A proof of this fact can be found in \cite{dispensedivillani} as 
long as we deal with the case $p>1$, while for $p=1$ we refer to \cite{ambpra}.

Let $B=t^{-1}(A)$. Up to modifying $t$ on the $\mu$-negligible set $A$, we may suppose $B\cap A=\emptyset$. 
Set $\mu_1=1_{B^c}\cdot\mu+1_{A\setminus N}\cdot\nu$; it is a probability measure with density $u_1$ given 
by $1_{B^c}u+1_{A}v=1_{B^c\setminus A}u+1_{A}v$ (this equality comes from $u=0$ on $A$). We have
$$F(\mu_1)=\int_{B^c\setminus A}f(u)\,d\lcal^n+\int_{A}f(v)\,d\lcal^n\leq F(\mu)+||f(v)||_{\infty}|\Omega|<+\infty.$$
Setting
$$t^*(x)=\begin{cases}t(x) &\text{if }x\in(A\cup B)^c,\\
                   x&\text{if }x\in(A\cup B),\end{cases}$$
we can see that $t^*$ is a transport map between $\mu_1$ and $\nu$. In fact, for 
any Borel set $E\subset\Omega$, we may express $(t^*)^{-1}(E)$ as the disjoint union 
of $E\cap A$, $E\cap B$, and $t^{-1}(E)\cap B^c\cap A^c$. Thus,
\begin{eqnarray*}
\mu_1((t^*)^{-1}(E))&=&\nu(E\cap A)+\nu(E\cap B\cap A)+\mu(t^{-1}(E)\cap B^c\cap A^c)\\
&=&\nu(E\cap A)+\mu(t^{-1}(E\cap A^c))=\nu(E),
\end{eqnarray*}
where we used the fact that $A\cap B=\emptyset$ and that $A^c$ is a set of full measure for $\mu$.
Consequently,
\begin{equation}\label{disstretta}
T_p(\mu_1,\nu)\leq\int_{(A\cup B)^c}\!\!\!|x-t(x)|^pu(x)dx<\int_{\Omega}\!|x-t(x)|^pu(x)dx=T_p(\mu,\nu).
\end{equation}
From this it follows that for $\mu_1$ Lemma \ref{stimafacileconconv} is not satisfied, 
since the integral term $\int_{\Omega}f'(u)(u_1-u)d\lcal^n$ is nonpositive, because $u_1>u$ 
only on $A$, where $f'(u)$ vanishes. The strict inequality in (\ref{disstretta}) follows from 
the fact that if $\int_{A\cup B}|x-t(x)|^pu(x)dx=0$, then for a.e. $x\in B$ it holds $u(x)=0$ or $x=t(x)$, 
which, by definition of $B$, implies $x\in A$; in both cases we are led to $u(x)=0$. This would give $\nu(A)=\mu(B)=0$, 
contradicting the assumptions $|A|>0$ and $v>0$ a.e. in $\Omega$.\qquad
\end{proof}

We need some results from duality theory in mass transportation that can be found 
in \cite{dispensedivillani}. In particular, we point out the notation of $c$-transform 
(a kind of generalization of the well-known Legendre transform): given a function $\chi$ on $\Omega$ 
we define its $c$-transform (or $c$-conjugate function) by
$$\chi^c(y)=\inf_{x\in\Omega}c(x,y)-\chi(x).$$
We will generally use $c(x,y)=|x-y|^p$.

\begin{teo}\label{peroraprincipale}
Under the same hypotheses of Lemma {\rm \ref{muriempietuttosenilofa}}, assuming also that $p>1$, 
if $\mu$ is optimal for $\mathfrak{F}^p_{\nu}$ and we denote by $\psi$ the unique, up to 
additive constants, Kantorovich potential for the transport between $\mu$ and $\nu$, there 
exists a constant $l$ such that the following relation holds:
\begin{equation}\label{reldenspot}
u=k(l-\psi)\; \text{a.e. in }\Omega .
\end{equation}
\end{teo}
\unskip

\begin{proof}
Let us choose an arbitrary measure $\mu_1$ such that $F(\mu_1)<+\infty$ (for instance a probability with bounded density $u_1$) and define $\mu_{\ve}=\mu+\ve(\mu_1-\mu)$. Let us denote by $\psi_{\ve}$ a Kantorovich potential between $\mu_{\ve}$ and $\nu$, chosen so that all the functions $\psi_{\ve}$ vanish at a same point. 
We can use the optimality of $\mu$ to write
$$T_p(\mu_{\ve},\nu)+F(\mu_{\ve})-T_p(\mu,\nu)-F(\mu)\geq 0.$$
By means of the duality formula, as $T_p(\mu_{\ve},\nu)=\int\psi_{\ve}d\mu_{\ve}+\int\psi_{\ve}^cd\nu$ and 
$T_p(\mu,\nu)\geq\int\psi_{\ve}d\mu+\int\psi_{\ve}^cd\nu$, we can write
$$\int\psi_{\ve}d(\mu_{\ve}-\mu)+F(\mu_{\ve})-F(\mu)\geq 0.$$
Recalling that $\mu_{\ve}-\mu=\ve(\mu_1-\mu)$ and that 
$$F(\mu_{\ve})-F(\mu)=\int \left(f(u+\ve(u_1-u))-f(u)\right)\,d\lcal^n,$$
we can divide by $\ve$ and pass to the limit. We know from Lemma \ref{lemmapsive} that $\psi_{\ve}$ 
converge towards the unique Kantorovich potential $\psi$ for the transport between $\mu$ and $\nu$. 
For the limit of the $F$ part we use Fatou's Lemma, as in Lemma \ref{stimafacileconconv}.
We then obtain at the limit
$$\int_{\Omega}(\psi(x)+f'(u(x)))(u_1(x)-u(x))\,dx\geq 0.$$
Notice that, since $u_1$ is arbitrary, one can infer that both $f'(u)$ and $uf'(u)$ are $L^1$ functions (to do so, take $u_1$ equal to $\max\{u/2, m\}$ for a suitable constant $m$ so that $\int u_1\,dx=1$: this gives $\int uf(u)\,dx<+\infty$).

If we now restrict ourselves to probabilities $\mu_1$ with bounded density $u_1$ we have 
\begin{equation*}\label{danumerare}
\int(\psi(x)+f'(u(x)))u_1(x)\,dx\geq\int(\psi(x)+f'(u(x)))u(x)\,dx.
\end{equation*}
Define first $l=\text{\rm{ess}}\inf_{x\in\Omega}\psi(x)+f'(u(x))$.
The left-hand side, by properly choosing $u_1$, can be made as close to $l$ as we want. 
Then we get that the function $\psi+f'(u)$, which is $\lcal^n$-a.e., and so also $\mu$-a.e., 
greater than $l$, integrated with respect to the probability $\mu$ gives a result less than or equal to $l$. It follows that
$$\psi(x)+f'(u(x))=l,\quad \mu{\mbox{\rm -a.e.}}\,x\in\Omega.$$
Together with the fact that, by Lemma \ref{muriempietuttosenilofa}, $u>0$ a.e., 
we get an equality valid $\lcal^n$-a.e., and so it holds that 
$$f'(u)=l-\psi.$$
We can then compose with $k$ and get the thesis.\qquad
\end{proof}

To establish Lemma \ref{lemmapsive}, which we used in the proof of Theorem \ref{peroraprincipale}, 
we have first to point out the following fact. In the transport between two probabilities, 
if we look at the cost $c(x,y)=|x-y|^p$ with $p>1$, there exists just one Kantorovich potential, 
up to additive constants, provided the absolutely continuous part of one of the measures has 
strictly positive density a.e. in the domain $\Omega$.
\begin{lem}\label{lemmapsive}
Let $\psi_{\ve}$ be Kantorovich potentials for the transport between $\mu_{\ve}=\mu+\ve(\mu_1-\mu)$ 
and $\nu$, all vanishing at a same point $x_0\in\Omega$. Suppose that $\mu=u\cdot\lcal^n$ and $u>0$ 
a.e. in $\Omega$, and let $\psi$ be the unique Kantorovich potential between $\mu$ and $\nu$ vanishing 
at the same point; then $\psi_{\ve}$ converge uniformly to $\psi$.
\end{lem}
\unskip

\begin{proof}
First, notice that the family $(\psi_{\ve})_{\ve}$ is equicontinuous since any $c$-concave function 
with respect to the cost $c(x,y)=|x-y|^p$ is $pD^{p-1}$-Lipschitz continuous (and Kantorovich potentials 
are optimal $c$-concave functions in the duality formula). Moreover, thanks to $\psi_{\ve}(x_0)=0 $, 
we also get equiboundedness and thus, by the Ascoli--Arzel\`{a} theorem, the existence of uniform limits up 
to subsequences. Let $\overline{\psi}$ be one of these limits, arising from a certain subsequence. 
From the optimality of $\psi_{\ve}$ in the duality formula for $\mu_{\ve}$ and $\nu$ we have, for 
any $c$-concave function $\ffi$,
$$\int\psi_{\ve}\,d\mu_{\ve}+\int\psi_{\ve}^c\,d\nu\geq\int\ffi\,d\mu_{\ve}+\int\ffi^c\,d\nu.$$
We want to pass to the limit as $\ve\tto 0 $: we have uniform convergence of $\psi_{\ve}$ but we need 
uniform convergence of $\psi_{\ve}^c$ as well. To get it, note that 
\begin{gather*}
\psi_{\ve}^{c}(x)=\inf_y |x-y|^p-\psi_{\ve}(y),\quad\overline{\psi}^{c}(x)=\inf_y |x-y|^p-\overline{\psi}(y),\\
|\psi_{\ve}^{c}(x)-\overline{\psi}^{c}(x)|\leq||\psi_{\ve}-\overline{\psi}||_{\infty}.
\end{gather*}
Passing to the limit as $\ve\tto 0$ along the considered subsequence we get, for any $\ffi$, 
$$\int\overline{\psi}\,d\mu+\int\overline{\psi}^c\,d\nu\geq\int\ffi\,d\mu+\int\ffi^c\,d\nu.$$
This means that
$\overline{\psi}$ is a Kantorovich potential for the transport between $\mu$ and $\nu$. 
Then, taking into account that $\overline{\psi}(x_0)=0$, we get the equality $\overline{\psi}=\psi$. 
Then we derive that the whole sequence converges to $\psi$.\qquad
\end{proof}

We now highlight that the relation we have proved in Theorem \ref{peroraprincipale} enables 
us to choose a density $u$ which is continuous. Moreover, it is also continuous in a quantified way, 
since it coincides with $k$ composed with a Lipschitz function with a fixed Lipschitz constant. 
As a next step we will try to extend such results to the case of general $\nu$ and then to the case $p=1$. 
The uniform continuity property we proved will be essential for an approximation process.

In order to go through our approximation approach, we need the following lemma, requiring the 
well-known theory of $\Gamma$-convergence. For all details about this theory, we refer to \cite{introgammaconve}.
\begin{lem}\label{gammanuep}
Given a sequence $(\nu_h)_h$ of probability measures on $\Omega$, supposing $\nu_h\deb\nu$ and $p>1$, 
it follows that the sequence of functionals $(\mathfrak{F}^p_{\nu_h})_h$ $\Gamma$-converges to 
the functional $\mathfrak{F}^p_{\nu}$ with respect to weak-$*$ topology on $\pical(\Omega)$. 
Moreover, if $\nu$ is fixed and we let $p$ vary, we have  $\Gamma$-convergence, according to the same topology, 
of the functionals $\mathfrak{F}^p_{\nu}$ to the functional $(\mathfrak{F}^{1}_{\nu})$ as $p\tto 1$.
\end{lem}
\unskip

\begin{proof}
For the first part of the statement, just notice that the Wasserstein distance is a metrization 
of weak-$*$ topology: consequently, since $T_p(\mu,\nu)=W_p^p(\mu,\nu)$, as $\nu_h\deb\nu$ we 
have uniform convergence of the continuous functionals $T_p(\cdot,\nu_h)$. This implies $\Gamma$-convergence 
and pointwise convergence. In view of Proposition 6.25 in \cite{introgammaconve}, concerning $\Gamma$-convergence 
of sums, we achieve the proof. The second assertion follows the same scheme
once we notice that, for each $p>1$ and every pair $(\mu,\nu)$ of probability measures, it holds that
$$W_1(\mu,\nu)\leq W_p(\mu,\nu)\leq D^{1-{\frac{1}{p}}}W_1^{{\frac{1}{p}}}(\mu,\nu).$$
This gives uniform convergence of the transport term, as
\begin{eqnarray*}
T_p(\mu,\nu)-T_1(\mu,\nu)&\leq& (D^{p-1}-1)T_1(\mu,\nu)\\
&\leq& D(D^{p-1}-1)\tto 0.\\
T_p(\mu,\nu)-T_1(\mu,\nu)&\geq& T_1^p(\mu,\nu)-T_1(\mu,\nu)\\
&\geq& (p-1)c(T_1(\mu,\nu))\geq \bar{c}\,(p-1)\tto 0,
\end{eqnarray*}
where $c(t)=t\log t$, $\bar{c}=\inf c$, and we used the fact $T_1(\mu,\nu)\leq D$.\qquad
\end{proof}

We now state in the form of lemmas two extensions of Theorem \ref{peroraprincipale} 
\begin{lem}\label{nugenericapgrande}
Suppose $p>1$ and fix an arbitrary $\nu\in\pical(\Omega)$. If $\mu$ is optimal for $\mathfrak{F}^p_{\nu}$, 
then there exists a Kantorovich potential $\psi$ for the transport between $\mu$ and $\nu$ such that 
{\rm \eqref{reldenspot}} holds.
\end{lem}
\unskip

\begin{proof}
We choose a sequence $(\nu_h)_h$ approximating $\nu$ in such a way that each $\nu_h$ satisfies 
the assumptions of Theorem \ref{peroraprincipale}. By Lemma \ref{gammanuep} and the properties 
of $\Gamma$-convergence, the space $\pical(\Omega)$ being compact and the functional $\mathfrak{F}^p_{\nu}$ 
having an unique minimizer (see, for instance, Chapter 7 in \cite{introgammaconve}), we get that $\mu_h\deb\mu$, 
where each $\mu_h$ is the unique minimizer of $\mathfrak{F}^p_{\nu_h}$. Each measure $\mu_h$ 
is absolutely continuous with density $u_h$. We use \eqref{reldenspot} to express $u_h$ in terms of 
Kantorovich potentials $\psi_h$ and get uniform continuity estimates on $u_h$. We would like to extract 
converging subsequences by the Ascoli--Arzel\`{a} theorem, but we also need equiboundedness. We may obtain 
this by using together the integral bound $\int u_h d\lcal^n=\int k(-\psi_h) d\lcal^n=1$ and the equicontinuity. 
So, up to subsequences, we have the following situation:
\begin{gather*}\mu_h=u_h\cdot\lcal^n,\qquad u_h=k(-\psi_h),\\
u_h\tto u,\qquad\psi_h\tto \psi  \text{ uniformly,}\\
\mu_h\deb\mu,\quad\mu=u\cdot\lcal^n,\qquad\nu_h\deb\nu,\end{gather*}
where we have absorbed the constants $l$ into the Kantorovich potentials. Clearly it is sufficient 
to prove that $\psi$ is a Kantorovich potential between $\mu$ and $\nu$ to reach our goal.

To see this, we consider that for any $c$-concave function $\ffi$, it holds that
$$\int\psi_{h}\,d\mu_h+\int\psi_h^c\,d\nu_h\geq\int\ffi\,d\mu_h+\int\ffi^c\,d\nu_h.$$
The thesis follows passing to the limit with respect to $h$, as in Lemma \ref{lemmapsive}.\qquad
\end{proof}

The next step is proving the same relation when $\nu$ is generic and $p=1$. We are in the same situation 
as before, and we simply need approximation results on Kantorovich potentials in the more difficult situation 
when the cost functions $c_p(x,y)=|x-y|^p$ vary with $p$.

\begin{lem}\label{nugenericap1}
Suppose $p=1$ and fix an arbitrary $\nu\in\pical(\Omega)$. If $\mu$ is optimal for $\mathfrak{F}^1_{\nu}$, 
then there exists a Kantorovich potential $\psi$ for the transport between $\mu$ and $\nu$ with cost $c(x,y)=|x-y|$ 
such that {\rm \eqref{reldenspot}} holds.
\end{lem}
\unskip

\begin{proof}
For any $p>1$ we consider the functional $\mathfrak{F}^p_{\nu}$ and its unique minimizer $\mu_p$. 
Thanks to Lemma \ref{nugenericapgrande} we get the existence of densities $u_p$ and Kantorovich potential 
$\psi_p$ between $\mu_p$ and $\nu$ with respect to the cost $c_p$, such that 
$$\mu_p=u_p\cdot\lcal^n,\quad u_p=k(-\psi_p).$$
By the Ascoli--Arzel\`a compactness result, as usual, we may suppose, up to subsequences, 
$$u_p\tto u,\quad\psi_p\tto\psi\text{  uniformly,}$$
and due to the $\Gamma$-convergence result in Lemma \ref{gammanuep}, since $\mathfrak{F}_{\nu}^1$ 
has an unique minimizer denoted by $\mu$ we also get
$$\mu_p\deb\mu,\qquad\mu=u\cdot\lcal^n.$$
As in Lemma \ref{nugenericapgrande}, we simply need to prove that $\psi$ is a Kantorovich potential between 
$\mu$ and $\nu$ for the cost $c_1$. The limit function $\psi$ is Lipschitz continuous with Lipschitz 
constant less than or equal to $\liminf_{p\tto 1}pD^{p-1}=1$, since it is approximated by $\psi_p$. 
Consequently $\psi$ is $c$-concave for $c=c_1$. We need to show that it is optimal in the duality formula.

Let us recall that, for any real function $\ffi$ and any cost function $c$, 
it holds that $\ffi^{cc}\geq\ffi$ and $\ffi^{cc}$ is a $c$-concave function whose $c$-transform is $\ffi^{ccc}=\ffi^c$. 
Consequently, by the optimality of $\psi_p$, we get
\begin{equation}\label{psipphicc}
\int\psi_pd\mu_p+\int\psi_p^{c_p}d\nu\geq\int\ffi^{c_p c_p}d\mu_p+\int\ffi^{c_p}d\nu
\geq\int\ffi d\mu_p+\int\ffi^{c_p}d\nu.
\end{equation}
We want to pass to the limit in the inequality between the first and the last term.
We start by proving that, for an arbitrary sequence $(\ffi_p)_p$, if $\ffi_p\tto\ffi_1$, we have the uniform convergence 
$\ffi_p^{c_p}\tto\ffi_1^{c_1}$. Let us take into account that we have uniform convergence on 
bounded sets of 
$c_p(x,y)=|x-y|^p$ to $c_1(x,y)=|x-y|$. Then we have
\begin{gather*}
\ffi_p^{c_p}(x)=\inf_y |x-y|^p-\ffi_p(y),\quad\ffi_1^{c_1}(x)=\inf_y |x-y|-\ffi_1(y),\\
|\ffi_p^{c,p}(x)-\ffi_1^{c,1}(x)|\leq||c_p-c_1||_{\infty}+||\ffi_p-\ffi_1||_{\infty},
\end{gather*}
which gives us the convergence we needed. We then obtain, passing to the limit as $p\tto 1$ in \eqref{psipphicc},
$$\int\psi\,d\mu+\int\psi^{c_1}\,d\nu\geq\int\ffi\,d\mu+\int\ffi^{c_1}\,d\nu.$$
By restricting this inequality to all $\ffi$ which are $c_1$-concave, we get that $\psi$ is a Kantorovich potential 
for the transport between $\mu$ and $\nu$ and the cost $c_1$.\qquad
\end{proof}

We can now state the main theorem of this section, whose proof consists only of putting 
together all the results we have obtained above.

\begin{teo}\label{formulafighissima}
Let $\Omega$ be a bounded domain in $\R^n$, $f$ be a $C^1$ strictly convex function, $p\geq 1$, and $\nu$ be
a probability measure on $\Omega$. Then there exists a unique measure $\mu\in\pical(\Omega)$ 
minimizing $\mathfrak{F}^p_{\nu}$ and it is absolutely continuous with density $u$. Moreover, 
there exists a Kantorovich potential $\psi$ for the transport between $\mu$ and $\nu$ and 
the cost $c(x,y)=|x-y|^p$ such that $u=k(-\psi)$ holds, where $k=(f')^{-1}$.
\end{teo}

Consequences on the regularity of $u$ come from this expression, which gives Lipschitz-type continuity, 
and from the relationship between Kantorovich potentials and optimal transport, which can be expressed 
through some PDEs. It is not difficult, for instance, in the case $p=2$, to obtain a Monge--Amp\`{e}re 
equation for the density $u$.

As we have already mentioned, we provide a sketch of an alternative proof to Theorem {\rm \ref{formulafighissima}} (a ``convex analysis proof'')
The idea of such a proof consists of looking at the subdifferential of the functional $\F^p_{\nu}$ in
order to get optimality conditions on the unique minimizer measure $\mu$ and its density $u$ (here we will
identify any absolutely continuous probability measure with its density).


{\it Sketch of Proof.~Step} 1. Consider the minimizing probability $\mu$ with density $u\in
L^1(\Omega)$ and define the vector space $X= span\left(L^{\infty}(\Omega),\left\{u\right\}\right)$, with
dual
$$X'=\left\{\xi\in L^1(\Omega):\,\int_{\Omega}|\xi|u\,d\lcal^n<+\infty\right\}.$$ 
Then we consider the
minimization problem for the functional $H$ defined on $X$ by
$$H(v)=\begin{cases}\F^p_{\nu}(v)&\text{if }v\in\pical(\Omega),\\
+\infty&\text{otherwise.}\end{cases}$$
It is clear that $u$ minimizes $H$. We will prove
\begin{equation}\label{sottoh}
\partial H(u)=\left\{f'(u)+\psi:\,\psi \text{ maximizes } \int_{\Omega}\phi d\mu +
\int_{\Omega}\phi^c d\nu \text{ for }\phi\in X'\right\}
\end{equation}
and then consider as an optimality condition $0\in\partial H(u)$. Notice that this works easily if one adds the assumption that $f$ has a polynomial growth, so that $f'(u)\in X'$. For the general case, as one may notice from the proof of Theorem \ref{peroraprincipale}, we can prove $f'(u)\in X'$ at least for the minimizer $u$. The subdifferential $\partial H$ of the
convex functional $H$ is to be considered in the sense of the duality between $X$ and $X'$. Notice that, 
in this setting, the $c$-transform $\phi^c$ of a function $\phi\in X'$ has to be defined replacing 
the $\inf$ with an $\text{\rm{ess}}\inf$. Finally, in
order to achieve the proof, it is sufficient to recognize that for a function
$\psi$ attaining the maximum in the duality formula, it necessarily holds that $\psi = \psi^{cc}$ a.e. on
$\left\{u>0\right\}$ and that this, together with $0=f'(u)+\psi$, implies $\psi=\psi^{cc}\wedge 0$. This
means that $\psi$ is an optimal $c$-concave function (since it is expressed as an infimum of two
$c$-concave functions) in the duality formula between $\mu$ and $\nu$, and so it is a Kantorovich
potential. In this way the thesis of Theorem \ref{formulafighissima} is achieved, provided
\eqref{sottoh} is proved.

{\it Step} 2. By using the same computations as in Lemma \ref{stimafacileconconv}, for
any $u_1\in X\cap\pical(\Omega)$, if we set $u_{\ve}=u+\ve(u_1-u)$, we may prove that
$$\lim_{\ve\tto 0}\frac{F(\mu_{\ve})-F(\mu)}{\ve}=\int_{\Omega}f'(u)(u_1-u)\,d\lcal^n.$$
Notice that, since $\int_{\Omega}f'(u)|u_1-u|\,d\lcal^n<+\infty$, by choosing $u_1=1/|\Omega|$ it follows
that $f'(u)$ and $f'(u)u$ are $L^1$ functions; i.e., $f'(u)\in X'$. Then it is possible to prove that this
implies $\partial H(u)=f'(u)+\partial T(u)$, where $T$ is the convex functional $T_p(\cdot,\nu)$.


{\it Step} 3. It remains to prove that
\begin{equation}\label{sottot}
\partial T(u)=\left\{\psi:\,\psi \text{ maximizes } \int_{\Omega}\phi d\mu +
\int_{\Omega}\phi^c d\nu\text{ for }\phi\in X'\right\}.
\end{equation}
In fact, if we define $K(\phi)=\int_{\Omega}\phi^c d\nu$, the key point is to prove that $K$ is concave
and upper semicontinuous in $\phi$. Then, by standard convex analysis tools, \eqref{sottot} is a
consequence of the equality $T(v)=\sup_{\phi}v\cdot\phi+K(\phi)$, where
$v\cdot\phi$ stands for the duality product between $X$ and $X'$ and equals $\int_{\Omega}v\phi\,d\lcal^n$.

\section{Applications to urban planning problems (with atomic services)}\label{sec4}

In this section we go through the consequences that Theorem \ref{formulafighissima} 
has in the problem of minimizing $\mathfrak{F}^p$, when this functional is built by using 
a term $G$ as in (\ref{defiG}), which forces the measure $\nu$, representing services, to be purely atomic. 
We have two goals: trying to have an explicit expression for $u$ in the case of a bounded domain $\Omega$ 
and proving an existence result in the case $\Omega=\R^n$.

\begin{teo}
Suppose $(\mu,\nu)$ is optimal for problem {\rm (\ref{e1minpb})}. Suppose also that the function $g$ is 
locally Lipschitz in $]0,1]$: then $\nu$ has finitely many atoms and is of the form $\nu=\sum_{i=1}^ma_i\delta_{x_i}$. 
\end{teo}
\unskip

\begin{proof}
It is clear that $\nu$ is purely atomic, i.e., a countable sum of Dirac masses. We want to show their finiteness. 
Consider $a=\max a_i$ (such a maximum exists since $\lim_i a_i=0$ and $a_i>0$) and let $L$ be the Lipschitz constant 
of $g$ on $[a,1]$. Now consider an atom with mass $a_i$ and modify $\nu$ by moving its mass onto the atom $x_j$ 
whose mass $a_j$ equals $a$, obtaining a new measure $\nu'$. The $G$-part of the functional decreases, while it 
may happen that the transport part increases. Since we do not change $\mu$, the $F$-part remains the same. 
By optimality of $\nu$ we get $T_p(\mu,\nu)+G(\nu)\leq T_p(\mu,\nu')+G(\nu')$ and thus
$$g(a_i)-La_i\leq g(a_i)+g(a)-g(a+a_i)\leq T_p(\mu,\nu')-T_p(\mu,\nu)\leq a_i D.$$
This implies
$$\frac{g(a_i)}{a_i}\leq D+L,$$
and by the assumption on the behavior of $g$ at $0$, this gives a lower bound $\delta$ on $a_i$. 
Since we have proved that every atom of $\nu$ has a mass greater than $\delta$, we may conclude that 
$\nu$ has finitely many atoms.\qquad
\end{proof}

Now we can use the results from last section.
\begin{teo}\label{formula quasi esplicita}
For any $\nu\in\pical(\Omega)$ such that $\nu$ is purely atomic and composed by finitely 
many atoms at the points $x_1,\dots,x_m$, if $\mu$ minimizes $\mathfrak{F}^p_{\nu}$, 
there exist constants $c_i$ such that
\begin{equation}\label{colmassimo}
u(x)=k\left((c_1-|x-x_1|^p)\vee\cdots\vee(c_m-|x-x_m|^p)\vee 0\right).
\end{equation}
In particular the support of $u$ is the intersection with $\Omega$ of a finite union 
of balls centered around the atoms of $\nu$.
\end{teo}
\unskip

\begin{proof}
On the Kantorovich potential $\psi$ appearing in Theorem \ref{formulafighissima}, we know that 
\begin{gather*}
\psi(x)+\psi^c(y)=|x-y|^p\quad\forall(x,y)\in spt(\gamma),\\
\psi(x)+\psi^c(y)\leq|x-y|^p\quad\forall(x,y)\in\Omega\times\Omega,
\end{gather*}
where $\gamma$ is an optimal transport plan between $\mu$ and $\nu$. Taking into account that $\nu$ 
is purely atomic we obtain, defining $c_i=\psi^c(x_i)$,
\begin{gather*}
-\psi(x)=c_i-|x-x_i|^p\quad\,\mu{\mbox{\rm -a.e.}}\, x\in\Omega_i,\\
-\psi(x)\geq c_i-|x-x_i|^p\quad\forall x\in\Omega,\,\forall i,
\end{gather*}
where $\Omega_i=t^{-1}(x_i)$, where $t$ is an optimal transport map between $\mu$ and $\nu$. Since $\mu-$a.e.~point  
in $\Omega$ is transported to a point $x_i$, we know that $u=0$ a.e. in the complement of $\bigcup_i\Omega_i$. 
Since, by $f'(u)=-\psi$, it holds that $-\psi(x)\geq 0$, one gets that everywhere in $\Omega$ the function $-\psi$ 
is greater than each of the terms $c_i-|x-x_i|^p$ and $0$, while a.e. it holds equality with at least one of them. 
By changing $u$ on a negligible set, one obtains \eqref{colmassimo}. The support of $\mu$, consequently, turns out 
to be composed of the union of the intersection with $\Omega$ of the balls $B_i=B(x_i,c_i^{1/p})$.\qquad
\end{proof}

Theorem \ref{formula quasi esplicita} allows us to have an almost explicit formula for the density of $\mu$. 
Formula \eqref{colmassimo} becomes more explicit when the balls $B_i$ are disjoint. We now give a sufficient 
condition on $\nu$ under which this fact occurs.

\begin{lem}\label{distanza fra gli atomi}
There exists a positive number $\overline{R}$, depending on the function $k$, such that any of the balls $B_i$ 
has a radius not exceeding $\overline{R}$. In particular, for any atomic probability $\nu$ such that the distance 
between any two of its atoms is larger than $2\overline{R}$, the balls $B_i$ are disjoint.
\end{lem}
\unskip

{\it Proof}.
Set $R_i=c_i^{1/p}$ and notice that
$$1=\int_{\Omega}u\geq\int_{B_i}k(c_i-|x-x_i|^p)\,dx=\int_0^{R_i}k(R_i^p-r^p)n\omega_n r^{n-1}\,dr,$$
where the number $\omega_n$ stands for the volume of the unit ball in $\R^n$. This inequality gives the 
required upper bound on $R_i$, since
$$\int_0^{R_i}k(R_i^p-r^p)n\omega_n r^{n-1}\,dr\geq C\int_0^{R_i-1}n r^{n-1}\,dr= C(R_i-1)^n.\qquad\endproof $$
When the balls $B_i$ are disjoint, we have $B_i=\Omega_i$ for every $i$ and we get a simple relation 
between radii and masses corresponding to each atom. The constants $c_i$ can then be found 
by using $R_i=c_i^{1/p}$. In fact, by imposing the equality of the mass of $\mu$ in the ball and of $\nu$ 
in the atom, the radius $R(m)$ corresponding to a mass $m$ satisfies
\begin{equation}\label{espliR}
m=\int_0^{R(m)}k(R(m)^p-r^p)n\omega_nr^{n-1}dr.
\end{equation}
For instance, if $f(s)=s^2/2$, we have
$$R(m)=\left(\frac{m(n+p)}{\omega_n p}\right)^{{\frac{1}{n+p}}}.$$
The second aim of this section is to obtain an existence result for the problem (\ref{e1minpb}) 
when $\Omega=\R^n$. A difference from the bounded case is the fact that we must look for 
minimization among all pairs of measures in $\mathcal{W}_p(\R^n)$, the $p$-th Wasserstein 
metric space (i.e., the space of measures $\lambda\in\pical(\R^n)$ such that $\int |x|^p\lambda(dx)<+\infty$, 
endowed with the distance $W_p$), rather than in $\pical(\R^n)$.

We start with some simple results about the minimization problem for $\mathfrak{F}^p_{\nu}$.
\begin{lem}\label{comportamento a nu fissato}
For every fixed $\nu\in\pical(\R^n)$ there exists a (unique if $f$ is strictly convex) 
minimizer $\mu$ for $\mathfrak{F}^p_{\nu}$: it belongs to $\mathcal{W}_p(\R^n)$ if and only if $\nu\in\mathcal{W}_p(\R^n)$, 
and if $\nu$ does not belong to this space, the functional $\mathfrak{F}^p_{\nu}$ 
is infinite on the whole $\mathcal{W}_p(\R^n)$. Moreover, if $\nu$ is compactly supported, the same happens for $\mu$.
\end{lem}
\unskip

\begin{proof}
The existence of $\mu$ comes from the direct method of the calculus of variations and 
the fact that if $(T_p(\mu_h,\nu))_h$ is bounded, then $(\mu_h)_h$ is tight. The  behavior 
of the functional with respect to the space $\mathcal{W}_p(\R^n)$ is trivial. Finally, the 
last assertion can be proved by contradiction, supposing $\mu(B(0,R)^c)>0$ for every $R<+\infty$ 
and replacing $\mu$ with $$\mu_R=1_{B_R}\cdot\mu+\frac{\mu(B_R^c)}{|B_r|}1_{B_r}\cdot\lcal^n,$$ 
where $B(0,r)$ is a ball containing the support of $\nu$. By optimality, we should have 
\begin{equation}\label{sptopt}
T_p(\mu_R,\nu)+F(\mu_R)\geq T_p(\mu,\nu)+F(\mu),
\end{equation}
but we have
\begin{gather}\label{disT_p}
T_p(\mu_R,\nu)-T_p(\mu,\nu)\leq -((R-r)^p-(2r)^p)\mu(B_R^c),\\
\label{disF}
F(\mu_R)-F(\mu)\leq \int_{B_r}\left[f\left(u+\frac{\mu(B_R^c)}{|B_r|}\right)-f(u)\right]\,d\lcal^n.
\end{gather}
By summing up \eqref{disT_p} and \eqref{disF}, dividing by $\mu(B_R^c)$, and taking
into account \eqref{sptopt}, we get
\begin{equation}\label{sommateedivise}
-((R-r)^p-(2r)^p)+ \frac{1}{\mu(B_R^c)}\int_{B_r}\left[f\left(u+\frac{\mu(B_R^c)}{|B_r|}\right)-f(u)\right]\,d\lcal^n\geq 0.
\end{equation}
Yet, by passing to the limit as $R\tto+\infty$ and $\mu(B_R^c)\tto 0$, the first
term in \eqref{sommateedivise} tends to $-\infty$, while the second is decreasing as
$R\tto +\infty$. This last one tends to $\int_{B_r}f'(u)d\lcal^n$, provided it is
finite for at least a value of $R$ (which ensures the finiteness of the limit as
well). To conclude, it is sufficient to prove that
$$\int_{B_r}\left[f\left(u+\frac{\mu(B_R^c)}{|B_r|}\right)-f(u)\right]\,d\lcal^n<+\infty.$$ 
This is quite easy in the case $f(z)=Az^q$ with $q>1$, while for general $f$ the assertion 
comes from the fact that $u$ is continuous on $\overline{B_r}$ and hence bounded. If $u=0$ a.e.~in $B_r$,
this is trivial; otherwise take the probability measures $\mu'=1_{B_r}/\mu(B_r)\cdot \mu$ and $\nu'=t_{\sharp}\mu'$ 
for an optimal transport map $t$ between $\mu$ and $\nu$. 
It is clear that $\mu'$ minimizes $\F^p_{\nu'}$ in the new domain $\Omega'=\overline{B_r}$. 
Then we may apply Theorem \ref{formulafighissima} and get the continuity of its density, 
which ensures the continuity of $u$ on $\overline{B_r}$.\qquad
\end{proof}

To go through our proof we need to manage minimizing sequences, in the sense of the following lemma.
\begin{lem}\label{successionespeciale}
It is possible to choose a minimizing sequence $((\mu_h,\nu_h))_h$ 
in $\mathcal{W}_p(\R^n)\times\mathcal{W}_p(\R^n)$ such that for every $h$ the measure $\nu_h$ 
is finitely supported, and the density of $\mu_h$ is given by {\rm (\ref{colmassimo})}, with disjoint 
balls centered at the atoms of $\nu_h$.
\end{lem}
\unskip

\begin{figure}[hbtp]
\begin{center}
\includegraphics[height=6.4cm]{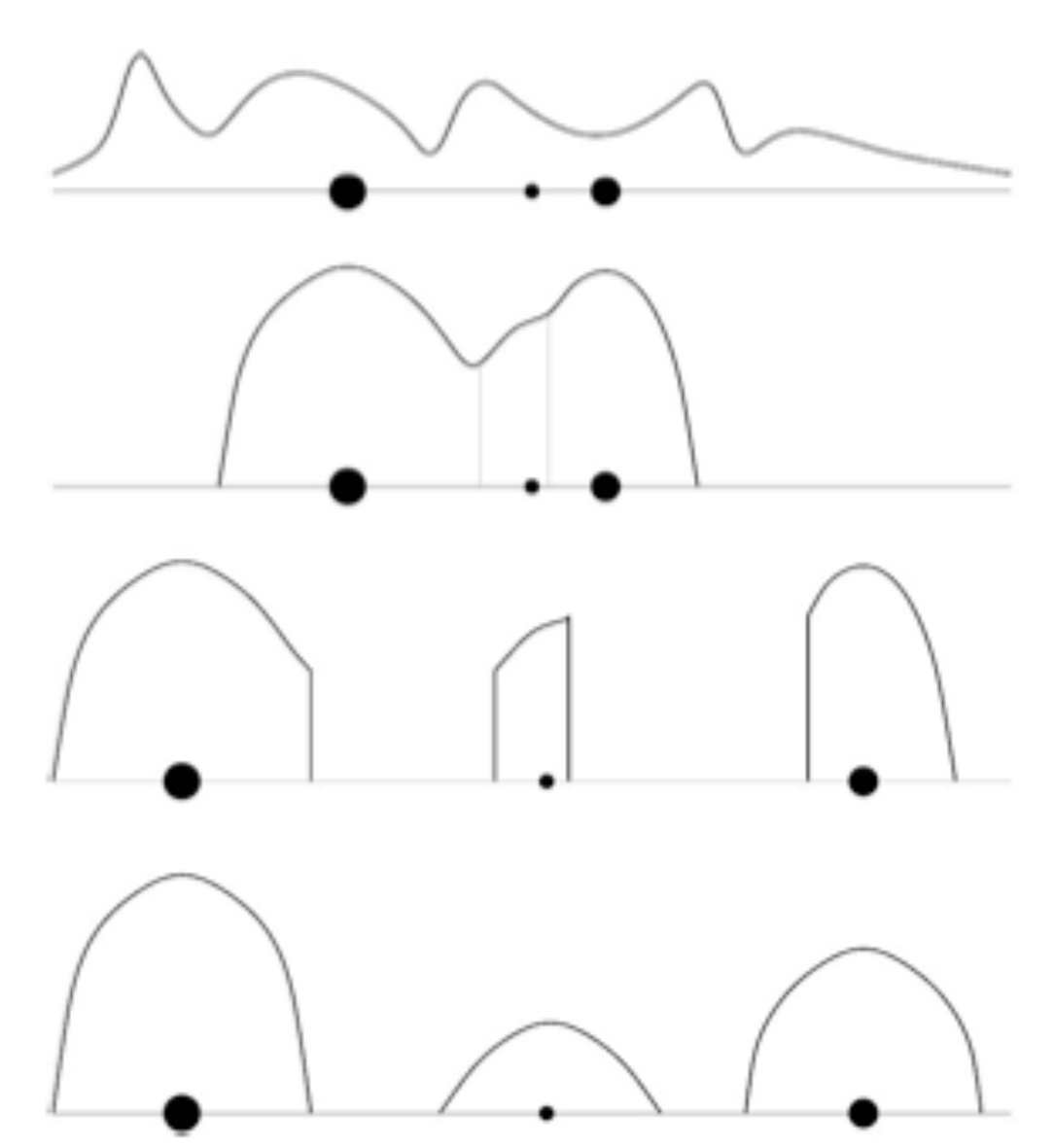}
\end{center}
\caption{Getting a well-behaved minimizing sequence.}
\label{improvements}
\end{figure}

\begin{proof}
First we start from an arbitrary minimizing sequence $((\mu'_h,\nu'_h))_h$. 
Then we approximate each $\nu'_h$ in $\mathcal{W}_p$ by a finite support 
measure $\nu''_h$. To do this we truncate the sequence of its atoms and 
move the mass in excess to the origin. In this way, we have 
$G(\nu''_h)\leq G(\nu'_h)$, by the subadditivity of $g$, while the value 
of the transport term increases by an arbitrary small quantity. 
Consequently, $((\mu'_h,\nu''_h))_h$ is still a minimizing sequence. 
Then we replace $\mu'_h$ by $\mu''_h$, chosen in such a way that it minimizes $\mathfrak{F}^p_{\nu''_h}$. 
By Lemma \ref{comportamento a nu fissato}, each $\mu''_h$ has a compact support. Then we translate 
every atom of each $\nu''_h$, together with its own set $\Omega_i$, to some disjoint sets $\Omega_i^*$. 
In this way we get new measures $\mu'''_h$ and $\nu'''_h$. The value of the functional in this step 
has not changed. We may choose to place the atoms of each $\nu'''_h$ so far from each other that 
each distance between atoms is at least $2\overline{R}$. Then we minimize again in $\mu$, getting a 
new sequence of pairs $((\mu''''_h,\nu'''_h))_h$, and we set $\nu_h=\nu'''_h$ and $\mu_h=\mu''''_h$. 
Thanks to Theorem \ref{formula quasi esplicita} and Lemma \ref{distanza fra gli atomi} the requirements 
of the thesis are fulfilled. The subsequent improvements of the minimizing sequence are shown in Figure \ref{improvements}.
\end{proof}

It is clear now that if one can obtain a uniform estimate on the number of atoms of the measures $\nu_h$, 
the existence problem is easily solved. In fact we already know that each ball belonging 
to the support of $\mu_h$ is centered at an atom of $\nu_h$ and has a radius not larger than $\overline{R}$. 
Provided we are able to prove an estimate like $\sharp\left\{\text{atoms of }\nu_h\right\}\leq N$, it would 
be sufficient to act by translation on the atoms and their corresponding balls, obtaining a new minimizing 
sequence (the value of $\F^p$ does not change) with supports all contained in a same bounded set (for instance, 
the ball $B_{N\overline{R}}$).

We now try to give sufficient conditions in order to find minimizing sequences where the number of atoms stays bounded. 
Notice that on sequences of the form given by Lemma \ref{successionespeciale}, the functional $\mathfrak{F}^p$ 
has the expression
\begin{equation}\label{espliF}
\mathfrak{F}^p(\mu_h,\nu_h)=\sum_{i=1}^{k(h)}E(m_{i,h}),\quad\text{if }\nu_h=\sum_{i=1}^{k(h)}m_{i,h}\delta_{x_{i,h}},
\end{equation}
where the quantity $E(m)$ is the total contribute given by an atom with mass $m$ to the functional. We may compute
\begin{multline}\label{espliE}
E(m)=g(m)+\!\int_0^{R(m)}\!\left[f(k(R(m)^p-r^p))+k(R(m)^p-r^p)r^p\right]n\omega_nr^{n-1}dr,
\end{multline}
taking into account the particular form of the density in the ball.

\begin{teo}\label{conunpocodizucchero}
Let us suppose $f\in
C^2((0,+\infty))$, and $g\in C^2((0,1])\cap C^0([0,1])$, in addition to all previous assumptions. 
Then the minimization problem for $\mathfrak{F}^p$ in $\mathcal{W}_p(\R^n)\times\mathcal{W}_p(\R^n)$ 
has a solution, provided
$$\limsup_{R\tto
0^+}g''\left(\int_0^{R}\!k(R^p-r^p)n\omega_nr^{n-1}dr\right)\int_0^R\!k'(R^p-r^p)n\omega_nr^{n-1}dr<-1.$$
\end{teo}
\unskip

\begin{proof}
According to what has been previously proven, it is sufficient to produce a minimizing 
sequence of the form of Lemma \ref{successionespeciale} with a bounded number of atoms. 
We claim that it is enough to prove that the function $E$ is subadditive on an interval $[0,m_0]$. 
In fact, having proven it, we start from a sequence $((\mu_h,\nu_h))_h$ built as in 
Lemma \ref{successionespeciale} and use the characterization of $\F^p$ given in \eqref{espliF}. 
Then we modify our sequence by replacing in each $\nu_h$ any pair of atoms of mass less than $m_0/2$ 
with a single atom with the sum of the masses. We keep atoms far away from each other in order 
to use \eqref{espliF}. We may perform such a replacement as far as we find more than one atom whose 
mass is less than or equal to $m_0/2$. At the end we get a new pair $((\mu'_h,\nu'_h))_h$, where the number 
of atoms of $\nu'_h$ is less than $N=1+\left\lfloor 2/m_0\right\rfloor$. The value of the functional $\F^p$ 
has not increased, thanks to the subadditivity of $E$ on $[0,m_0]$. 

Taking into account that $E(0)=0$ and that concave functions vanishing at $0$ are subadditive, 
we look at concavity properties of the function $E$ in an interval $[0,m_0]$. It is sufficient 
to compute the second derivative of $E$ and find it negative in a neighborhood of the origin. 

By means of the explicit formula (\ref{espliE}) and also taking into account (\ref{espliR}), 
setting $E(m)=g(m)+K(R(m))$, we start by computing $dK/dr$. Using the facts that $f'\,\circ\,k=id$ and that $k(0)=0$, 
we can obtain the formula
$$\frac{dK(R(m))}{dm}(m)=R(m)^p.$$
From another derivation and some standard computation we finally obtain
$$E''(m)=g''(m)+\frac{1}{\int_0^{R(m)}k'(R(m)^p-r^p)n\omega_nr^{n-1}dr}.$$
The assumption of this theorem ensures that such a quantity is negative for small $m$, and so the proof is achieved.\qquad
\end{proof}

\begin{oss}{\rm Notice that when the functions $f$ and $g$ are of the form $f(t)=at^q,\,q>1,\,g(t)=bt^r,\,r<1$, 
with $a$ and $b$ positive constants, it holds that
\begin{gather*}
g''\left(\int_0^{R}\!k(R^p-r^p)n\omega_nr^{n-1}dr\right)\leq
-C R^{(n+\frac{p}{q-1})(r-2)},\\
\int_0^R\!k'(R^p-r^p)n\omega_nr^{n-1}dr\leq
CR^{n+p\frac{2-q}{q-1}}, 
\end{gather*}
and so the $\limsup $ in Theorem \ref{conunpocodizucchero} may be estimated from above by
$$\lim_{R\tto 0^+}-C R^{\frac{p}{q-1}(r-q)+n(r-1)}=-\infty.$$
Consequently the assumption in Theorem \ref{conunpocodizucchero} is always verified when $f$ and $g$ are power functions.}
\end{oss}

\begin{oss}{\rm 
From the proof of the existence theorem it is clear that there exists a minimizing pair 
$(\mu,\nu)\in\mathcal{W}_p(\R^n)\times\mathcal{W}_p(\R^n)$ where $\nu$ has finitely many 
atoms and $\mu$ is supported in a finite, disjoint union of balls centered at the atoms of 
$\nu$ and contained in a bounded domain $\Omega_0$, with a density given by Theorem \ref{formula quasi esplicita}. 
The same happens if we look for the minimizers in a bounded domain $\Omega$, provided $\Omega$ is large enough 
to contain $\Omega_0$, and hence a solution to the problem in $\R^n$. For example, all the open sets containing 
$N$ balls of radius $\overline{R}$ admit a minimizing solution supported in disjoint balls.}
\end{oss}

\begin{figure}[h!]
   \begin{minipage}[b]{0.40\linewidth}
      \centering \includegraphics[height=3.6cm]{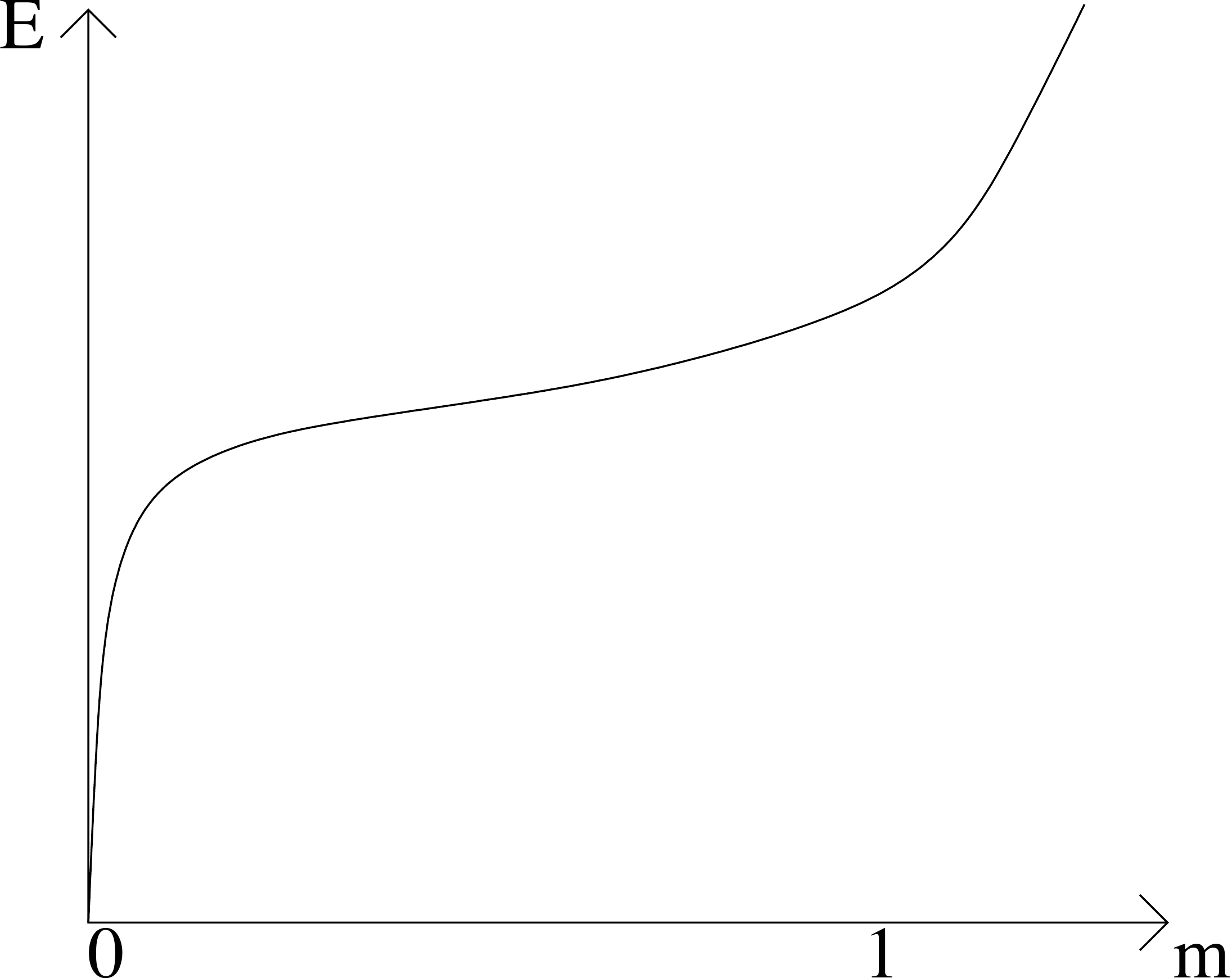}
\caption{Typical behavior of E.}
\label{grafico}
   \end{minipage}\hfill
   \begin{minipage}[b]{0.48\linewidth}   
      \centering\includegraphics[height=3.6cm]{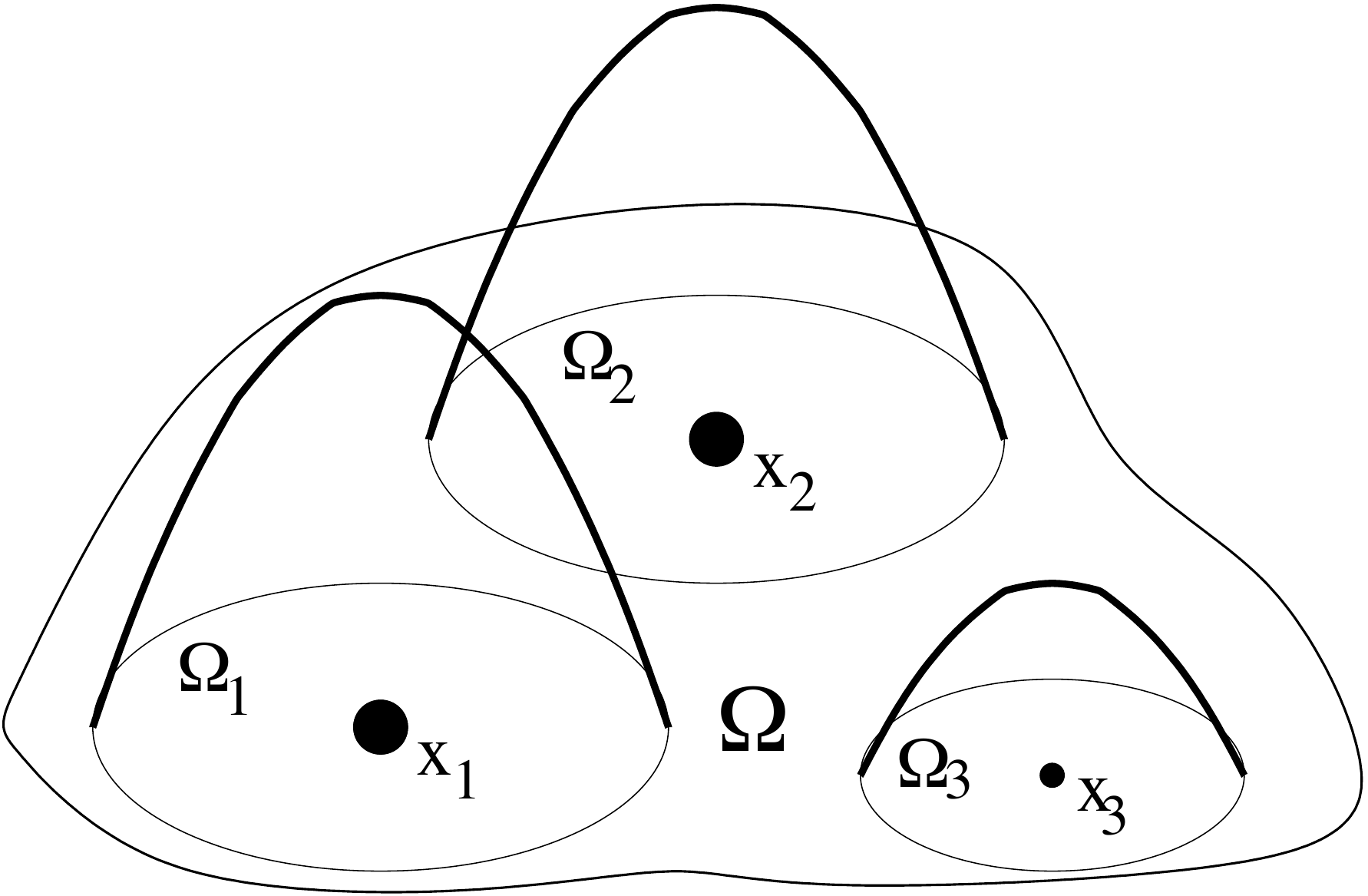}
\caption{Sketch of the solution in $R^n$.}
\label{solutionRn}

   \end{minipage}
\end{figure}

We conclude by stressing the fact that in order to solve the problem in $\R^n$, we have only 
to look at the function $E$ and find out the number of atoms and their respective masses 
$(m_i)_{i=1\dots k}$. The problem to solve is then
\begin{equation}\label{problema con E}
\min\left\{\sum_{i=1}^k E(m_i)\,:\,k\in\mathbb{N},\,\sum_{i=1}^k m_i=1\right\}.
\end{equation}
Typically, for instance when $f$ and $g$ are power functions, the function $E$ involved in \eqref{problema con E} is a concave-convex function, as sketched in Figure \ref{grafico}. Due to such a concave-convex behavior, in general it is not clear whether the values of the numbers $m_i$ solving \eqref{problema con E} and representing subcities' sizes are all equal or may be different (see also Figure \ref{solutionRn}).

\section*{Acknowledgments}
The authors wish to thank an anonymous referee for pointing out the alternative proof sketched in Section \ref{sec3}.


\end{document}